\documentclass[12pt]{amsart}

\theoremstyle{plain} %%%%%%%%%%%%%%%%%%%%%%%%%%%%%%%%%
\newtheorem{Theorem}{Theorem}[section]

\newtheorem{Corollary}[Theorem]{Corollary}

\theoremstyle{definition} %%%%%%%%%%%%%%%%%%%%%%%%%%%%
\newtheorem{Definition}[Theorem]{Definition}

\numberwithin{equation}{section}

\newcommand{\alg}[1]{\mathcal{#1}}                % algebras

\newcommand{\abs}[1]{\left\lvert #1 \right\rvert}
\newcommand{\intpart}{\alg{I}}  % connected partitions
\newcommand{\Piconn}{\Pi^{conn}}  % connected partitions
\newcommand{\Piconnpair}{\Pi^{conn,pair}}  % connected partitions
\newcommand{\Piirr}{\Pi^{irr}}  % irreducible partitions
\newcommand{\NC}{NC}  % connected partitions
\newcommand{\NCirr}{\NC^{irr}}  % irreducible partitions
\renewcommand{\phi}{\varphi}
\begin{document}

\title{Free Cumulants and Enumeration of Connected Partitions}
\author{Franz Lehner}

\address{
Franz Lehner\\
In\-sti\-tut f\"ur Mathe\-ma\-tik C\\
Tech\-ni\-sche Uni\-ver\-si\-t\"at Graz\\
Stey\-rer\-gas\-se 30, A-8010 Graz\\
Austria}
%http://finanz.math.tu-graz.ac.at/~lehner
\email{lehner@finanz.math.tu-graz.ac.at}

\date{\today}

\begin{abstract}
  A combinatorial formula is derived which expresses
  free cumulants in terms of classical comulants.
  As a corollary, we give a combinatorial interpretation of free cumulants
  of classical distributions, notably Gaussian and Poisson
  distributions. The latter count connected pairings
  and connected partitions respectively.
  The proof relies on M\"obius inversion on the partition lattice.
\end{abstract}

\maketitle

\section{Introduction}
\label{sec:intro}

A partition of $[n]=\{1,2,\dots,n\}$ is called \emph{connected}
if no proper subinterval of $[n]$ is a union of blocks.
Connected partitions have been studied by various authors
under various names.
They were introduced as
\emph{irreducible diagrams} in \cite{Kleitman:1970:proportions}
and \cite{Hsieh:1973:proportions} (see also earlier work of Touchard
\cite{Touchard:1952:configuration}),
where their asymptotic enumeration properties are studied.
They are reconsidered under the name of \emph{linked diagrams} in
\cite{Stein:1978:LinkedDiagramsI},
\cite{SteinEverett:1978:LinkedDiagramsII},
and \cite{NijenhuisWilf:1979:enumeration}
and a recursion is derived there.
They also appear as a basic example in the theory of decomposable combinatorial
objects developed in \cite{Beissinger:1985:enumeration}.
There and in \cite{BenderRichmond:1984:asymptotic, BenderOdlyzkoRichmond:1985:asymptotic}
asymptotics of more general irreducible partitions are studied.

We reserve the term \emph{irreducible partition} to partitions which
cannot be ``factored'' into subpartitions, i.e.\ partitions of
$[n]$ for which $1$ and $n$ are in the same connected component.
A partition is called \emph{noncrossing} if its blocks do not intersect
in their graphical representation, i.e., if there are no two
distinct blocks $B_1$ and $B_2$ and elements $a, c\in B_1$ and $b,d\in B_2$
s.t.\ $a<b<c<d$. Equivalently one could say that a partition is noncrossing
if each of its connected components consists of exactly one block.
Typical examples of these types of partitions are shown in Fig.~\ref{fig:Partitions}.

We denote the lattice of partitions of $[n]$ by $\Pi_n$, the irreducible
partitions by $\Piirr_n$ and the order ideal of
connected partitions by $\Piconn_n$; the lattice of noncrossing partitions
will be denoted by $\NC_n$, and the sublattice of irreducible noncrossing partitions by
$\NCirr_n$.
Finally, let us denote by $\intpart_n$ the lattice of interval partitions,
i.e.\ the lattice of partitions consisting entirely of intervals.

%
% connected partition
%
\begin{center}
  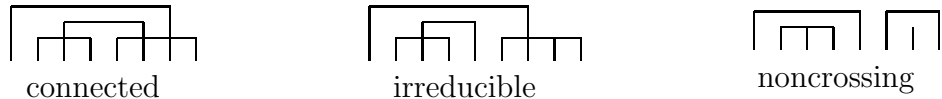
\begin{figure}[htbp]
    \begin{minipage}{.3\textwidth}
      \begin{center}
      % fprint(Unquoted,Text,"p.tex",SetPart::Cycle2TeX([[1, 7], [2, 4], [3, 6], [5, 8]], [10,6]))
        \begin{picture}(80,20.4)(1,0)% [[1, 7], [2, 4], [3, 6], [5, 8]]
          \put(10,0){\line(0,1){20.4}}
          \put(20,0){\line(0,1){8.4}}
          \put(30,0){\line(0,1){14.4}}
          \put(40,0){\line(0,1){8.4}}
          \put(50,0){\line(0,1){8.4}}
          \put(60,0){\line(0,1){14.4}}
          \put(70,0){\line(0,1){20.4}}
          \put(80,0){\line(0,1){8.4}}
          \put(20,8.4){\line(1,0){20}}
          \put(30,14.4){\line(1,0){30}}
          \put(10,20.4){\line(1,0){60}}
          \put(50,8.4){\line(1,0){30}}
        \end{picture}
        
        connected
      \end{center}
    \end{minipage}
%
% irreducible partition
%
    \begin{minipage}{.3\textwidth}
      \begin{center}
      % fprint(Unquoted,Text,"p.tex",SetPart::Cycle2TeX([[1, 7], [2, 4], [3, 5], [6, 8, 9]], [10,6]))
      \begin{picture}(90,20.4)(1,0)% [[1, 7], [2, 4], [3, 5], [6, 8, 9]]
        \put(10,0){\line(0,1){20.4}}
        \put(20,0){\line(0,1){8.4}}
        \put(30,0){\line(0,1){14.4}}
        \put(40,0){\line(0,1){8.4}}
        \put(50,0){\line(0,1){14.4}}
        \put(60,0){\line(0,1){8.4}}
        \put(70,0){\line(0,1){20.4}}
        \put(80,0){\line(0,1){8.4}}
        \put(90,0){\line(0,1){8.4}}
        \put(20,8.4){\line(1,0){20}}
        \put(30,14.4){\line(1,0){20}}
        \put(10,20.4){\line(1,0){60}}
        \put(60,8.4){\line(1,0){30}}
      \end{picture}

        irreducible
      \end{center}
    \end{minipage}
%
% noncrossing partition
%
    \begin{minipage}{.3\textwidth}
      % fprint(Unquoted,Text,"p.tex",SetPart::Cycle2TeX([[1, 5], [2, 3, 4], [6, 8],[7]], [10,6]))
      \begin{center}
        \begin{picture}(80,14.4)(1,0)% [[1, 5], [2, 3, 4], [6, 8], [7]]
          \put(10,0){\line(0,1){14.4}}
          \put(20,0){\line(0,1){8.4}}
          \put(30,0){\line(0,1){8.4}}
          \put(40,0){\line(0,1){8.4}}
          \put(50,0){\line(0,1){14.4}}
          \put(60,0){\line(0,1){14.4}}
          \put(70,0){\line(0,1){8.4}}
          \put(80,0){\line(0,1){14.4}}
          \put(20,8.4){\line(1,0){20}}
          \put(10,14.4){\line(1,0){40}}
          \put(70,8.4){\line(1,0){0}}
          \put(60,14.4){\line(1,0){20}}
        \end{picture}
        
        noncrossing
      \end{center}
    \end{minipage}
    \label{fig:Partitions}
    \caption{Typical partitions}
  \end{figure}
\end{center}

\section{Incidence algebras}
\label{sec:incidencealgebras}

Before recalling more facts about partitions, let us briefly introduce
the main concepts about posets and incidence algebras which will be needed
in the sequel.
Rota et al.\  \cite{DoubiletRotaStanley:1972:foundationsVI}
introduced the \emph{reduced incidence algebra} of a poset.
Let $(P,\leq)$ be a finite poset. On the space $I(P)$
of complex-valued functions $f(x,y)$
defined on the pairs $(x,y)$ s.t.\ $x\leq y$
(``triangular matrices'') we introduce a convolution 
(``multiplication of triangular matrices'') by
\begin{equation}
  \label{eq:IncidenceConvolution}
  f*g(x,y) = \sum_{x\leq z\leq y} f(x,z)\,g(z,y)  
\end{equation}
With this operations $I(P)$ becomes a unital algebra, the \emph{incidence algebra}
of the poset $P$, with identity 
$$
\delta(x,y) = 
\begin{cases}
  1 & \text{if $x=y$}\\
  0 & \text{if $x\ne y$}
\end{cases}
.
$$
It is clear by induction that a function is invertible under this convolution
if and only if all the ``diagonal'' entries $f(\pi,\pi)$ are nonzero.
One prominent invertible function is the zeta-function defined by
$$
\zeta(x,y) = 1 \quad \text{for $x\leq y$}
$$
Its inverse is called the \emph{M\"obius function} and satisfies the recursion
$$
\mu(x,y) =
\begin{cases}
  1                            & \text{if $x=y$}\\
  -\sum_{x\leq z < y} \mu(x,z) & \text{if $x<y$}
\end{cases}
$$
Then we have $\zeta*\mu=\mu*\zeta=\delta$ and more generally 
the \emph{M\"obius inversion formula} holds:
for any pair of functions $f(x)$ and $g(x)$ on $P$ we have the following equivalences 
\begin{align}
  \label{eq:MoebiusinversionEquivalence:ge}
  \forall x: f(x) = \sum_{y\geq x} g(y) 
  &\iff
  \forall x: g(x) = \sum_{y\geq x} \mu(x,y)\,f(y) \\
  \label{eq:MoebiusinversionEquivalence:le}
  \forall x: f(x) = \sum_{y\leq x} g(y)
  &\iff
  \forall x: g(x) = \sum_{y\leq x} f(y)\,\mu(y,x)
\end{align}

In our case the posets $P$ will be the partition lattice $\Pi_n$ and some of its sublattices,
namely the lattice $\intpart_n$ of interval partitions and the lattice $\NC_n$ of
non-crossing partitions.
Denote by $L_n$ any one of these lattices,
then every segment $[\pi,\sigma]$ is canonically isomorphic to a finite
direct product of full lattices $L_1^{k_1}\times L_2^{k_2}\times \cdots$.
The sequence of exponents $(k_1,k_2,\dots)$ will be called the \emph{type} of the segment
$[\pi,\sigma]$.
The \emph{reduced incidence algebra} is the algebra of functions whose values
on an interval $[\pi,\sigma]$ only depend on the type of the interval.
An even smaller class of functions is the set of \emph{multiplicative functions}
whose values at $[\pi,\sigma] \simeq L_1^{k_1}\times L_2^{k_2}\times \cdots$
are given by
$$
f(\pi,\sigma) = f_1^{k_1} f_2^{k_2} \cdots
$$
where $(f_j)$ is a given sequence of numbers.
These functions are defined on all partition lattices $\Pi_n$ simultaneously.
Examples of such functions are the zeta function and the M\"obius function and
it can be shown that the space of multiplicative functions is closed under the
convolution \eqref{eq:IncidenceConvolution}.

Multiplicative functions on the lattice
of noncrossing partitions were studied by Speicher in \cite{Speicher:1994:multiplicative}
and interval partitions in \cite{SpeicherWoroudi:1997:boolean},
see also \cite{vonWaldenfels:1973:approach} and \cite{vonWaldenfels:1975:interval}.

We refer to 
\cite{Aigner:1979:combinatorial} or
\cite[ch.~3]{Stanley:1986:enumerative}
for more information on incidence algebras.

\section{Cumulants}

Cumulants linearize convolution of probability measures coming from
various notions of independence.
\begin{Definition}
  A \emph{non-commutative probability space} is pair $(\alg{A},\phi)$
  of a (complex) unital algebra $\alg{A}$ and a unital linear functional $\phi$.
  The elements of $\alg{A}$ are called \emph{(non-commutative) random variables}.
  The collection of moments $\mu_n(a) = \phi(a^n)$
  of such a random variable $a\in\alg{A}$ will be called
  its \emph{distribution} and denoted $\mu_a = (\mu_n(a))_n$.
\end{Definition}
Thus noncommutative probability theory follows the general ``quantum'' philosophy
of replacing function algebras by noncommutative algebras.
We will review several notions of independence below.
Convolution is defined as follows. Let $a$ and $b$ be ``independent'' random variables.
Then the convolution of the distributions of $a$ and $b$ is defined to be
the distribution of the sum $a+b$. In all the examples below, the distribution
of the sum of  ``independent'' random variables only depends on the individual
distributions of the summands and therefore convolution is well defined and
the $n$-th moment $\mu_n(a+b)$ will be a polynomial function of the moments
of $a$ and $b$ of order less or equal to $n$.

For our purposes it is sufficient to axiomatize cumulants as follows.
\begin{Definition}
  \label{def:cumulants}
  Given a notion of independence on a noncommutative probability space $(\alg{A},\phi)$,
  a sequence of maps $a\mapsto k_n(a)$, $n=1,2,\dots$ is called a \emph{cumulant sequence} if 
  it satisfies
  \begin{enumerate}
   \item $k_n(a)$ is a polynomial in the first $n$ moments of $a$ with leading term $\mu_n(a)$.
    This ensures that conversely the moments can be recovered from the cumulants.
   \item homogenity: $k_n(\lambda a) = \lambda^n k_n(a)$.
   \item additivity: if $a$ and $b$ are ``independent'' random variables,
    then $k_n(a+b)=k_n(a)+k_n(b)$.
  \end{enumerate}
\end{Definition}

M\"obius inversion on the lattice of partitions plays a crucial role
in the combinatorial approach to cumulants.
We will need three kinds cumulants here, corresponding to classical, free and boolean
independence, and which are connected to the three lattices
of partitions considered in section~\ref{sec:intro}.
Let $X$ be a random variable with distribution $\psi$
and moments $m_n = m_n(X) = \int x^n\, d\psi(x)$

\subsection{Classical cumulants}
  Let 
  $$
  \alg{F}(z) = \int e^{xz}\,d\psi(x)
       = \sum_{n=0}^\infty \frac{m_n}{n!}\,z^n
  $$
  be the formal Laplace transform
  (or exponential moment generating function).
  Taking the formal logarithm we can write this series as
  $$
  \alg{F}(z) = e^{K(z)}
  $$
  where 
  $$
  K(z) = \sum_{n=1}^\infty \frac{\kappa_n}{n!}\,z^n
  $$
  is the \emph{cumulant generating function} and the numbers
  $\kappa_n$ are called the \emph{(classical) cumulants}
  of the random variable $X$.

  Set partitions come in as follows.
  Let $f$ and $g$ be the multiplicative functions in the reduced incidence algebra
  of $\Pi_n$ determined by the sequence $m_n$ and $\kappa_n$ respectively,
  then $f=g*\zeta$ and $g=f*\mu$, i.e.,
  if for a partition $\pi=\{\pi_1,\pi_2,\dots,\pi_p\}$
  we put $m_\pi = m_{\abs{\pi_1}}m_{\abs{\pi_2}}\cdots m_{\abs{\pi_p}}$
  and
  $\kappa_\pi = \kappa_{\abs{\pi_1}}\kappa_{\abs{\pi_2}}\cdots \kappa_{\abs{\pi_p}}$,
  then we can express moments and cumulants mutually as
  $$
  m_\pi = \sum_{\sigma\leq\pi} \kappa_\sigma
  \qquad
  \kappa_\pi = \sum_{\sigma\leq\pi} m_\sigma \, \mu(\sigma,\pi)
  $$

  For example, the standard Gaussian distribution $\gamma=N(0,1)$
  has cumulants
  $$
  \kappa_n(\gamma)
  = \begin{cases}
      1 & n=2\\
      0 & n\ne 2
    \end{cases}
  $$
  while the Poisson distribution $P_\lambda$
  (with weights $P_\lambda(\{k\}) = e^{-\lambda} \frac{\lambda^k}{k!}$)
  has cumulants $\kappa_n(P_\lambda) = \lambda$.
  
  It follows that the even moments $m_{2n} = \frac{2n!}{2^n n!}$
  of the standard gaussian distribution
  count the number of pairings of a set with the corresponding number of elements.

  On the other hand, the moments of a Poisson variable with parameter $1$ are known
  as \emph{Bell numbers} $B_n$ and they are equal to the numbers of partitions of the finite
  sets with the corresponding cardinalities.
  The moment interpretation leads to \emph{Dobinski's formula}
  (cf.\ \cite{Rota:1964:number})
  $$
  B_n = e^{-1} \sum_{k=0}^\infty \frac{k^n}{k!}
  $$

\subsection{Free Cumulants}
  Free cumulants were introduced by Speicher
  \cite{Speicher:1994:multiplicative} in his combinatorial
  approach to Voiculescu's free probability theory
  \cite{VoiculescuDykemaNica:1992:free}.
  Given our random variable $X$, let
  \begin{equation}
    \label{eq:ordinarymgf}
    M(z) = 1 + \sum_{n=1}^\infty m_n z^n
  \end{equation}
  be its ordinary moment generating function.
  Define a formal power series 
  $$
  C(z) = 1 + \sum_{n=1}^\infty c_n z^n
  $$
  implicitly by the equation
  $$
  C(z) = C(zM(z))
  .
  $$
  Then the coefficients $c_n$ are called the \emph{free} or \emph{non-crossing cumulants}.
  The latter name stems from the fact that
  combinatorially these cumulants are obtained by M\"obius inversion on the lattice of
  non-crossing partitions:
  \begin{equation}
    \label{eq:NCcumulantsMoebiusinversion}
    m_\pi = \sum_{\substack{\sigma\in\NC_n\\ \sigma\leq\pi}} c_\sigma
    \qquad
    c_\pi = \sum_{\substack{\sigma\in\NC_n\\ \sigma\leq\pi}} m_\sigma \,\mu_{\NC}(\sigma,\pi)
  \end{equation}

\subsection{Boolean cumulants}
Boolean cumulants linearize \emph{boolean convolution} 
\cite{SpeicherWoroudi:1997:boolean,vonWaldenfels:1973:approach,vonWaldenfels:1975:interval}.
Let again $M(z)$ be the ordinary moment generating function of a random variable $X$
defined by \eqref{eq:ordinarymgf}.
It can be written as
$$
M(z) = \frac{1}{1-H(z)}
$$
where
$$
H(z) = \sum_{n=1}^\infty h_n z^n
$$
and the coefficients are called \emph{boolean cumulants}.
Combinatorially the connection between moments and boolean cumulants is described
by M\"obius inversion on the lattice of interval partitions:
\begin{equation}
  \label{eq:booleancumulantsMoebiusinversion}
  m_\pi = \sum_{\substack{\sigma\in\intpart_n\\ \sigma\leq\pi}}
           h_\sigma
  \qquad
  h_\pi = \sum_{\substack{\sigma\in\intpart_n\\ \sigma\leq\pi}}
           m_\sigma \,\mu_{\intpart}(\sigma,\pi)
\end{equation}
The term ``boolean cumulants'' is due to the fact that the lattice of interval partitions
is anti-isomorphic to the boolean lattice of subsets of the same set with the first
element removed. The isomorphism maps a partition to the set of first elements of its
blocks, where clearly the first element of the first block is always the same and therefore
redundant.

\section{Enumeration of connected partitions}

In this note we apply M\"obius inversion to show that the free cumulants
count connected partitions with certain weights given by the classical
cumulants. The result is inspired by \cite{Beissinger:1985:enumeration}.

\begin{Theorem}
  Let $(m_n)$ be a (formal) moment sequence with classical cumulants $\kappa_n$.
  Then the free cumulants of $m_n$ are equal to
  \begin{equation}
    \label{eq:freeintermsofclassical}
    c_n = \sum_{\pi\in \Piconn_n} \kappa_\pi    
  \end{equation}
  the boolean cumulants are equal to
  $$
  h_n = \sum_{\pi\in \Piirr_n} \kappa_\pi = \sum_{\pi\in \NCirr_n} c_\pi
  $$
\end{Theorem}
\begin{proof}
  We consider only the identity \eqref{eq:freeintermsofclassical},
  the proof of the others being similar
  (and also contained in the lattice path picture of \cite{Lehner:2001:LatticePaths}).
  For $\sigma\in\Pi_n$ we denote by $\bar{\sigma}$ its \emph{noncrossing closure},
  that is the smallest noncrossing partition $\pi$ s.t.\ $\sigma\leq\pi$.
  This noncrossing partition is obtained from $\sigma$ by taking unions of
  all blocks which cross in the graphical representation like in the example
  depicted in Fig.~\ref{fig:piandbarpi}.

  \begin{center}
    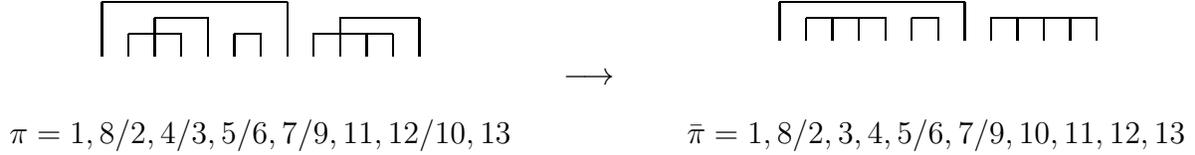
\begin{figure}[htbp]
      \begin{minipage}{.4\textwidth}
        \begin{center}
        % fprint(Unquoted,Text,"p.tex",SetPart::Cycle2TeX([[1, 8], [2, 4], [3, 5], [6, 7], [9, 11, 12], [10, 13]], [10,6]))
          \begin{picture}(130,20.4)(1,0)% [[1, 8], [2, 4], [3, 5], [6, 7], [9, 11, 12], [10, 13]]
            \put(10,0){\line(0,1){20.4}}
            \put(20,0){\line(0,1){8.4}}
            \put(30,0){\line(0,1){14.4}}
            \put(40,0){\line(0,1){8.4}}
            \put(50,0){\line(0,1){14.4}}
            \put(60,0){\line(0,1){8.4}}
            \put(70,0){\line(0,1){8.4}}
            \put(80,0){\line(0,1){20.4}}
            \put(90,0){\line(0,1){8.4}}
            \put(100,0){\line(0,1){14.4}}
            \put(110,0){\line(0,1){8.4}}
            \put(120,0){\line(0,1){8.4}}
            \put(130,0){\line(0,1){14.4}}
            \put(20,8.4){\line(1,0){20}}
            \put(30,14.4){\line(1,0){20}}
            \put(60,8.4){\line(1,0){10}}
            \put(10,20.4){\line(1,0){70}}
            \put(90,8.4){\line(1,0){30}}
            \put(100,14.4){\line(1,0){30}}
          \end{picture}
        \end{center}
        
        $$\pi =  1, 8/2, 4/3, 5/6, 7/9, 11, 12/10, 13$$
%        $$\pi =  \{\{1, 8\}, \{2, 4\}, \{3, 5\}, \{6, 7\}, \{9, 11, 12\}, \{10, 13\}\}$$
      \end{minipage}
      \hskip2em $\longrightarrow$ \hskip2em
      \begin{minipage}{.4\textwidth}
        % NCP::Cycle2TeX([[1,8],[2,3,4,5],[6,7],[9,10,11,12,13]],[10,6])
        \begin{center}
        \begin{picture}(130,14.4)(1,0)% [[1, 8], [2, 3, 4, 5], [6, 7], [9, 10, 11, 12, 13]]
          \put(10,0){\line(0,1){14.4}}
          \put(20,0){\line(0,1){8.4}}
          \put(30,0){\line(0,1){8.4}}
          \put(40,0){\line(0,1){8.4}}
          \put(50,0){\line(0,1){8.4}}
          \put(60,0){\line(0,1){8.4}}
          \put(70,0){\line(0,1){8.4}}
          \put(80,0){\line(0,1){14.4}}
          \put(90,0){\line(0,1){8.4}}
          \put(100,0){\line(0,1){8.4}}
          \put(110,0){\line(0,1){8.4}}
          \put(120,0){\line(0,1){8.4}}
          \put(130,0){\line(0,1){8.4}}
          \put(10,14.4){\line(1,0){70}}
          \put(20,8.4){\line(1,0){30}}
          \put(60,8.4){\line(1,0){10}}
          \put(90,8.4){\line(1,0){40}}
        \end{picture}
        \end{center}
        
        $$\bar{\pi} = 1,8 / 2,3,4,5 / 6,7 / 9,10,11,12,13$$
      \end{minipage}
      \caption{A partition $\pi$ and its noncrossing closure}
      \label{fig:piandbarpi}
    \end{figure}
  \end{center}

  For each $\pi\in\NC_n$ define the number
  $$
  \tilde{c}_\pi = \sum_{\substack{\sigma\in\Pi_n\\ \bar{\sigma}=\pi}} \kappa_\sigma
  $$
  Now note that the preimage of $\hat{1}_n = \{\{1,2,\dots,n\}\}$
  is the set of all connected partitions, i.e.\
  $$
  \tilde{c}_n := \tilde{c}_{\hat{1}_n}
  = \sum_{\sigma\in\Piconn_n} \kappa_\sigma
  $$
  
  For general $\pi\in\NC_n$ ,
  by considering subpartitions induced by the blocks of $\pi$,
  we have multiplicativity
  $\tilde{c}_\pi = \prod_{B\in\pi} \tilde{c}_{|B|}$.
  Now for $\pi\in\NC_n$ we can collect terms as follows
  \begin{align*}
    m_\pi &= \sum_{\substack{\rho\in \Pi_n\\ \rho\leq\pi}}
              \kappa_\rho \\
        &= \sum_{\substack{\sigma\in \NC_n\\ \sigma\leq\pi}}
            \sum_{\substack{\rho\in\Pi_n\\ \bar{\rho}=\sigma}}
             \kappa_\rho \\
        &= \sum_{\substack{\sigma\in \NC_n\\ \sigma\leq\pi}}
            \tilde{c}_\sigma
  \end{align*}
  and by M\"obius inversion \eqref{eq:MoebiusinversionEquivalence:le}
  and \eqref{eq:NCcumulantsMoebiusinversion}
  it follows that $c_\pi = \tilde{c}_\pi$.
\end{proof}

\begin{Corollary}
  The free cumulants of the standard gaussian variable are equal to
  the number of connected pairings.
  $$
  c_{2n}(\gamma) = \# \Piconnpair_{2n}
  $$
\end{Corollary}

\begin{Corollary}
  The free cumulants of the Poisson distribution with parameter $1$ are equal
  to the number of connected partitions
  $$
  c_{n}(P_1) = \# \Piconn_{n}
  .
  $$
  Moreover, if we leave the formal parameter $\lambda$, the expression
  for the $n$-th cumulants is the generating function of the numbers of blocks
  of the connected partitions.
  $$
  c_{n}(P_\lambda) = \sum_{\pi\in \Piconn_{n}} \lambda^{\abs{\pi}}
  $$
\end{Corollary}
Similar identities hold for the $q$-Gaussian and $q$-Poisson laws
\cite{BozejkoKummererSpeicher:1997:qGaussian,Nica:1995:qCumulants}
where the free cumulants provide a generating function
of the number of \emph{left-reduced crossings} of the connected partitions.
Alternatively, the free cumulants of the $q$-Poisson law of 
\cite{Anshelevich"2001:qCumulants,SaitohYoshida:2000:qPoissonOP,SaitohYoshida:2000:qPoissonFock}
count the number of \emph{reduced crossints}, cf.\ \cite{Biane:1997:properties}.
In all these examples a continued fraction expansion of the moment generating
function is known and the free cumulants can be expressed via Lagrange inversion.
See also \cite{Lehner:2001:LatticePaths}.

\providecommand{\bysame}{\leavevmode\hbox to3em{\hrulefill}\thinspace}

\end{document}